\documentclass[journal]{IEEEtran}
\usepackage{epsfig}

\usepackage{amssymb}

\usepackage{amsthm,diagbox}
\usepackage{graphicx}
\usepackage{amsfonts}
\usepackage{amssymb,bm}
\usepackage{amsmath}
\usepackage{subfigure}
\usepackage{float}
\usepackage{color}

\begin{document}

\title{A Smooth Double Proximal Primal-Dual Algorithm for a Class of Distributed Nonsmooth Optimization Problem}

\author{Yue~Wei,
		Hao~Fang,~\IEEEmembership{Member,~IEEE,}
		Xianlin~Zeng,~\IEEEmembership{Member,~IEEE,}
        Jie~Chen,~\IEEEmembership{Senior Member,~IEEE,}
        and~Panos M.~Pardalos
\thanks{This work was supported by Projects of Major International (Regional) Joint Research Program NSFC (Grant no. 61720106011), NSFC (Grant no. 61621063, 6157306261673058), Beijing Advanced Innovation Center for Intelligent Robots and Systems (Beijing Institute of Technology), Key Laboratory of Biomimetic Robots and Systems (Beijing Institute of Technology), Ministry of Education, Beijing, 100081, China. P. M. Pardalos was supported by the Paul and Heidi Brown Preeminent Professorship at ISE, University of Florida.}
\thanks{Y. Wei, H. Fang, X. Zeng and J. Chen are with the School of Automation, Beijing Institute of Technology, and also the Key Laboratory of Intelligent Control and Decision of Complex Systems, 100081, Beijing, China; P. M. Pardalos is with the Center of Applied Optimization (CAO), Industrial and Systems Engineering, University of Florida. (Corresponding author: H. Fang, email: fangh@bit.edu.cn) }}

\maketitle

\begin{abstract}
This technical note studies a class of distributed nonsmooth convex consensus optimization problem. The cost function is a summation of local cost functions which are convex but nonsmooth. Each of the local cost functions consists of a twice differentiable (smooth) convex function and two lower semi-continuous (nonsmooth) convex functions. We call this problem as \emph{single-smooth plus double-nonsmooth} (SSDN) problem. Under mild conditions, we propose a distributed double proximal primal-dual optimization algorithm. The double proximal operator is designed to deal with the difficulty caused by the unproximable property of the summation of those two nonsmooth functions. Besides, it can also guarantee that the proposed algorithm is locally Lipschitz continuous. An auxiliary variable in the double proximal operator is introduced to estimate the subgradient of the second nonsmooth function. Theoretically, we conduct the convergence analysis by employing Lyapunov stability theory. It shows that the proposed algorithm can make the states achieve consensus at the optimal point. In the end, nontrivial simulations are presented and the results demonstrate the effectiveness of the proposed algorithm.
\end{abstract}


\IEEEpeerreviewmaketitle

\section{Introduction}
This note aims to propose a smooth algorithm for a class of distributed nonsmooth convex consensus optimization problem. The cost function in this problem is a summation of local cost functions and each of them consists of a smooth (differentiable) convex function and two nonsmooth (nondifferentiable) convex functions, all of which are derived from practical meanings. Although both of the nonsmooth functions are proximable, their summation might not be (a function being proximable means the proximal operator of this function has a closed or semi-closed form solution and is computationally easy to evaluate $\cite{CTOS3}$). Moreover, all the solutions to this optimization problem must achieve consensus. We call such a problem as \emph{single-smooth plus double-nonsmooth} (SSDN) problem. A wide range of problems in the field of machine learning and multi-agent system optimization fall into SSDN form. For example, in the distributed version of the fused LASSO problem $\cite{Lasso}$, the least squares loss is smooth. The $l_{1}$ penalty and another penalty with respect to the changes of the temporal/spatial structure in this problem are usually nonsmooth. Besides, a consensual solution to this problem should be obtained. In the distributed constrained consensus optimization problems of multi-agent systems, which are different from $\cite{DCOR1}$-$\cite{Yipeng2}$, the local cost function of agent $i$ contains a smooth function and two nonsmooth functions standing for two different constraints. Consensus needs to be achieved at the optimal point of this problem. Due to the important implementation mentioned above, the SSDN problem has attracted increasing attentions.

There are two important categories of the existing algorithms for solving the distributed nonsmooth optimization problem. The first category includes the nonsmooth algorithms (see, $\cite{DCO1}$-$\cite{DCO5}$ for example.) whose convergence was proven based on nonsmooth analysis $\cite{NCO}$. Many distributed subgradient algorithms were presented as discrete-time systems with diminishing step-sizes $\cite{DCO1}$-$\cite{DCO3}$ or continuous-time systems $\cite{Xie1}$ to deal with the distributed nonsmooth optimization problems. $\cite{Yang1}$-$\cite{DCO5}$ proposed various projected  primal-dual algorithms in the form of differential inclusions to solve the distributed nonsmooth constrained optimization problems. Since the algorithms in the first category are nonsmooth, the subgradient of the nonsmooth cost function is discontinuous and may cause violent vibration in the systems, which is actually unacceptable especially for multi-agent dynamical systems. Moreover, it is also hard to prove the convergence properties of nonsmooth algorithms. The second category includes smooth algorithms (see, $\cite{DCFO1}$-$\cite{DCFO3}$ for example.) which employed the proximal methods $\cite{PM}$ to solve the distributed nonsmooth optimization problem. $\cite{DCFO1}$-$\cite{DCFO4}$ proposed several proximal gradient algorithms for solving the distributed optimization problems with nonsmooth local cost functions. $\cite{DCFO2}$-$\cite{DCFO3}$ designed proximal gradient alternating direction method of multipliers (PG-ADMM) to solve the distributed nonsmooth convex optimization problems. However, the aforementioned algorithms can not deal with SSDN problems (pointed out in $\cite{TOS1}$-$\cite{TOS3}$), since the summation of two proximable nonsmooth functions may not be proximable. As for the centralized version of SSDN problem, $\cite{CTOS3}$ proposed a double proximal gradient algorithm, but it is infeasible for the distributed nonsmooth optimization problem. 

In this note, the SSDN problem is addressed. A distributed smooth double proximal primal-dual algorithm is proposed to make the states of the agents achieve consensus at the optimal point to the SSDN problem. The contributions of this note are summarized as follows. 

\textbf{(i)} This note explores the SSDN problem, which widely arises in the subject of machine learning and multi-agent system optimization. To the best of our knowledge, very few studies on smooth algorithms for this problem have been carried out under the framework of distributed optimization.

\textbf{(ii)} A distributed \textbf{smooth} primal-dual algorithm using the double proximal operator is proposed. To tackle the main difficulty caused by the unproximable property of the summation of the two nonsmooth functions in each local cost function and ensure smoothness of the proposed algorithm, the double proximal operator is employed. In the double proximal operator, an auxiliary variable is designed to estimate the subgradient of the second nonsmooth function. Moreover, compared to the double proximal gradient algorithm $\cite{CTOS3}$, the proposed algorithm is a \textbf{distributed} primal-dual algorithm.

\textbf{(iii)} The convergence and correctness of the proposed algorithm are proved by using Lyapunov stability theory. The proof avoids using the nonsmooth analysis and provides novel insights into analysis of primal-dual type algorithms. 

The rest of the note is organized as follows. In Section II, some basic definitions of graph theory, and basic concepts of proximal operator are presented. Section III shows the SSDN problem. In Section IV, we propose a distributed smooth continuous-time primal-dual algorithm, which applies the double proximal operator. In Section V, the proof for the convergence and correctness of the algorithm is presented. In Section VI, simulations show the effectiveness of our proposed algorithm. Finally, Section VII concludes this note.

\section{Mathematical Preliminaries}
In this section, we introduce necessary notations, definitions and preliminaries about graph theory and proximal operator. 


\subsection{Graph Theory}
A weighted undirected graph $\mathcal{G}$ is denoted by $\mathcal{G(V,E,A)}$, where $\mathcal{V} = \lbrace 1, \dots, n \rbrace$ is a set of nodes, $\mathcal{E} = \lbrace (i, j) : i, j \in \mathcal{V}; i \neq j \rbrace \subset \mathcal{V} \times \mathcal{V}$ is a set of edges, and $\mathcal{A} = [a_{i,j}] \in \mathbb{R}^{n \times n}$ is a weighted adjacency matrix such that $a_{i,j} = a_{j,i} > 0$ if $(j,i) \in E$ and $a_{i,j} = 0$ otherwise, where $\mathbb{R}^{n \times n}$ denotes the set of $n$-by-$n$ real matrices. $j \in \mathcal{N}_{i}$ denote agent $j$ is a neighbour of agent $i$. The Laplacian matrix is $L_{n} = D - \mathcal{A}$, where $D \in \mathbb{R}^{n \times n}$ is diagonal with $D_{i,i} = \sum^{n}_{j=1} a_{i,j}$, $i \in \lbrace 1, \dots, n \rbrace$. Specifically, if the weighted graph $\mathcal{G}$ is undirected and connected, then $L_{n} = L^{T}_{n} \geq 0$, $rank$ $(L_{n}) = n - 1$ and $ker(L_{n}) = \lbrace k\textbf{1}_{n} : k \in \mathbb{R}\rbrace$, where $\mathbb{R}^{n}$ denotes the set of $n$-dimensional real column vectors and $\textbf{1}_{n} \in \mathbb{R}^{n}$ is the vector of all ones.

\subsection{Proximal Operator}
Let $f(\delta)$ be a lower semi-continuous convex function for $\delta \in \mathbb{R}^{r}$. Then the proximal operator $prox_{f}[\eta]$ and the Moreau envelope $M_{f}[\eta]$ of $f(\delta)$ at $\eta \in \mathbb{R}^{r}$ are
\begin{equation}
prox_{f}[\eta] = \arg \min_{\delta} \lbrace f(\delta) + \frac{1}{2} \Vert \delta - \eta \Vert^{2} \rbrace
\end{equation}

\begin{equation}
M_{f}[\eta] =  \inf_{\delta} \lbrace f(\delta) + \frac{1}{2} \Vert \delta - \eta \Vert^{2} \rbrace 
\end{equation}
where $\Vert \cdot \Vert$ denotes the Euclidean norm.

The Moreau envelope $M_{f}[\eta]$ is essentially a smooth or regularized form of $f(\delta)$ at $\eta$: it is continuously differentiable, even when $f(\delta)$ is not.

Define the indicator function of a closed convex set $\Omega$ as $I_{\Omega}(\delta) = 0$ if $\delta \in \Omega$  and $I_{\Omega}(\delta) = +\infty$ otherwise. We have $prox_{I_{\Omega}}[\eta] = P_{\Omega}[\eta]$, where $P_{\Omega}[\eta] = \arg \min_{\delta \in \Omega} \Vert \delta - \eta \Vert$ is the projection operator. Let $ \partial f(\delta)$ denote the subgradient of $f(\delta)$. If $f(\delta)$ is convex, then $\partial f(\delta)$ is monotone, that is, $(\zeta_{\delta_{1}} - \zeta_{\delta_{2}})^{T}(\delta_{1} - \delta_{2}) \geq 0$ for all $\delta_{1} \in \mathbb{R}^{r}, \delta_{2} \in \mathbb{R}^{r}$, $\zeta_{\delta_{1}} \in \partial f(\delta_{1})$, and $\zeta_{\delta_{2}} \in \partial f(\delta_{2})$. $\delta =
prox_{f} [\eta]$ is equivalent to
\begin{equation}
\eta - \delta \in \partial f(\delta)
\label{Proximal Property}
\end{equation}

\section{Problem Description}
In this section, the SSDN problem is formulated. We consider a network of $n$ agents with first-order dynamics, interacting over a graph $\mathcal{G}$. For each agent $i \in \lbrace 1, \cdots, n \rbrace$, there are three functions $f_{i}^{0} ,f_{i}^{1}, f_{i}^{2}: \mathbb{R}^{q} \to \mathbb{R}$ contained in the local cost function $f_{i}(x_{i}): \mathbb{R}^{q} \to \mathbb{R}$, where $f^{0}_{i}$ is a smooth convex function, $f^{1}_{i}(f^{2}_{i})$ is a nonsmooth convex function, and $\mathbb{R}$ denotes the set of real numbers. Each agent $i$ only has the information about $f_{i}^{0}$, $f_{i}^{1}$ and $f_{i}^{2}$. 

The SSDN problem is described as follows
\begin{eqnarray}
\min_{x \in \mathbb{R}^{nq}} F(x) \!\!\!\!& = \!\!\!\!& \sum_{i=1}^{n} f_{i}(x_{i}) \notag \\
{\rm where} \ f_{i}(x_{i}) \!\!\!\!& = \!\!\!\!& f^{0}_{i}(x_{i})+f^{1}_{i}(x_{i})+f^{2}_{i}(x_{i}) \notag\\
s.t. \quad L_{nq} x \!\!\!\!& = \!\!\!\!& \textbf{0}_{nq} \label{Problem 2}
\end{eqnarray}
where $L_{nq}= L_{n} \otimes I_{q}$, $I_{q}$ denotes the $q \times q$ identity matrix and $L_{n} \otimes I_{q}$ is the Kronecker product of matrices $L_{n}$ and $I_{q}$. $x_{i} \in \mathbb{R}^{q}$ is the state of $i$-th agent, and $x = [ x_{1}^{T}, x_{2}^{T}, \cdots, x_{n}^{T} ]^{T}$. The constraint presented in  ($\ref{Problem 2}$)  indicates that all the solutions must achieve consensus. Each agent only exchanges information with its neighbours in a fully distributed manner.

\emph{Remark 3.1:}
Problem ($\ref{Problem 2}$) is a very general model, which provides a new viewpoint of many problems in recent distributed constrained optimization $\cite{DCO1,DCO3,DCO2}$. For example, if $f^{1}_{i}(x_{i})$ is an indicator function of the convex set $\Omega_{i}$, then $f^{1}_{i}(x_{i})$ is equivalent to the constraint set $x_{i} \in \Omega_{i}$. If,  additionally, $f^{2}_{i}(x_{i}) = \mu \Vert x_{i} \Vert_{1}$, the optimization problem is a fused LASSO problem with the set constraints $x \in \cap_{i=1}^{n}\Omega_{i}$, where $\Vert \cdot \Vert_{1}$ denotes the $l_{1}$ norm.

The assumptions below are made for the wellposedness of the problem ($\ref{Problem 2}$) in this section.

\emph{Assumption 3.1:}  $f^{0}_{i}$ is twice continuously differentiable and strongly convex for all $i \in \lbrace 1, \cdots, n \rbrace$, which means there exists a constant $c > 0$ such that for agent $i$,
\begin{equation}
(\nabla f^{0}_{i}(\vartheta_{1}) - \nabla f^{0}_{i}(\vartheta_{2}))^{T}(\vartheta_{1} - \vartheta_{2}) \geq c \Vert \vartheta_{1} - \vartheta_{2} \Vert^{2}
\label{Strongly Convex}
\end{equation} 
where $\vartheta_{1} \in \mathbb{R}^{q}$, $\vartheta_{2} \in \mathbb{R}^{q}$, $\vartheta_{1} \neq \vartheta_{2}$, and $(\cdot)^{T}$ denotes transpose. Without loss of generality, we assume $c > 1$.

\emph{Remark 3.2:} If $0 < c \leq 1$, there exists a function $f^{0'}_{i}(x) = K f^{0}_{i}(x)$ for agent $i$ with $ K >\frac{1}{c}$ such that
\begin{equation}
(\nabla \! f^{0'}_{i}(\vartheta_{1}) - \nabla \! f^{0'}_{i}(\vartheta_{2}))^{T}\!(\vartheta_{1} - \vartheta_{2})\! \geq Kc \Vert \vartheta_{1} - \vartheta_{2} \Vert^{2} \! > \!\Vert \vartheta_{1} - \vartheta_{2} \Vert^{2}
\end{equation} 

\emph{Assumption 3.2:}  $f^{1}_{i}$ and $f^{2}_{i}$ are (nonsmooth) lower semi-continuous closed proper convex functions for all $i \in \lbrace 1, \cdots, n \rbrace$, and they are proximable.

\emph{Assumption 3.3:}  The weighted graph $\mathcal{G}$ is connected and undirected.

\emph{Assumption 3.4:}  There exists at least one finite solution to problem ($\ref{Problem 2}$).

Then, we arrive at the following lemma by the Karush-Kuhn-Tucker (KKT) condition of convex optimization problems.

\emph{Lemma 3.1:}
Under the Assumptions 3.1-3.4, a feasible point $x^{*} \in \mathbb{R}^{nq}$ is a minimizer to problem ($\ref{Problem 2}$) if and only if there exist $x^{*} = \textbf{1}_{n} \otimes w^{*} \in \mathbb{R}^{nq}$, $w^{*} \in \mathbb{R}^{q}$, and $v^{*} \in \mathbb{R}^{nq}$ such that
\begin{eqnarray}
& & \textbf{0}_{nq} \in \nabla F^{0}(x^{*}) + \partial F^{1}(x^{*}) + \partial F^{2}(x^{*}) + L_{nq} v^{*} \notag \\
& & L_{nq} x^{*} = \textbf{0}_{nq} \label{KKT}
\end{eqnarray}
where $\textbf{0}_{nq} \in \mathbb{R}^{nq}$ is the vector of all zeros, $\nabla F^{0}(x) = [(\nabla f^{0}_{1}(x_{1}))^{T}, (\nabla f^{0}_{2}(x_{2}))^{T}, \cdots, (\nabla f^{0}_{n}(x_{n}))^{T} ]^{T}$, $\partial F^{j}(x) = [(\partial f^{j}_{1}(x_{1}))^{T}, (\partial f^{j}_{2}(x_{2}))^{T}, \cdots, (\partial f^{j}_{n}(x_{n}))^{T} ]^{T}$, $j \in \lbrace 1, 2 \rbrace$. 
The proof of Lemma 3.1 is omitted since it is a trivial extension of the proof for Theorem 3.25 in $\cite{NO}$, 

\section{Algorithm Design}
In this section, we present a distributed smooth double proximal primal-dual algorithm for solving the problem ($\ref{Problem 2}$). In order to deal with the difficulty caused by the unproximable property of $f^{1}_{i}(x_{i})+f^{2}_{i}(x_{i})$ for each agent $i$, here we introduce an auxiliary variable $z(t) \in \mathbb{R}^{nq}$ combined with a constant parameter $\gamma \in \mathbb{R}^{+}$ to estimate $\partial F^{2}(x)$, such that there exists a feasible point $z^{*}$ satisfying that
\begin{eqnarray}
- \nabla F^{0}(x^{*}) - L_{nq} v^{*} + \gamma z^{*} \!\!\!\!& \in \!\!\!\!& \partial F^{1}(x^{*}) \notag \\
- \gamma z^{*} \!\!\!\!& \in \!\!\!\!& \partial F^{2}(x^{*}) \label{z}
\end{eqnarray}
where $\mathbb{R}^{+}$ denotes the set of positive real numbers, $x^{*}$ and $v^{*}$ are defined as same as in ($\ref{KKT}$).

According to ($\ref{KKT}$) and ($\ref{z}$), we propose a smooth algorithm as follows
\begin{eqnarray}
\dot{x}_{i}(t) \!\!\!\!& = \!\!\!\!& prox_{f^{1}_{i}} \Big[ x_{i}(t) \!-\! \nabla f^{0}_{i}(x_{i}(t)) \!-\! \alpha \! \sum_{j \in \mathcal{N}_{i}} \! a_{i,j}(v_{i}(t) - v_{j}(t))  \notag \\
\!\!\!\!& \!\!\!\!& - \alpha \sum_{j \in \mathcal{N}_{i}} a_{i,j} (x_{i}(t) - x_{j}(t)) + \gamma z_{i}(t) \Big] - x_{i}(t) \notag \\
\dot{z}_{i}(t) \!\!\!\!& = \!\!\!\!& prox_{f^{2}_{i}}[x_{i}(t) - \gamma z_{i}(t)] - x_{i}(t) \notag \\
\dot{v}_{i}(t) \!\!\!\!& = \!\!\!\!& \alpha \sum_{j \in \mathcal{N}_{i}} a_{i,j} (x_{i}(t) - x_{j}(t)) 
\label{Algorithm 1}
\end{eqnarray}
where $t \geq 0$, $i \in \lbrace 1, \cdots, n \rbrace$,  $v_{i}$ is the Lagrange multiplier for agent $i$, $0 < \alpha < \frac{1}{\lambda_{max}(L_{n})}$, $0 < \gamma < 1 - \alpha \lambda_{max}(L_{n})$, and $\lambda_{max}(L_{n})$ denotes the largest eigenvalue of Laplacian matrix $L_{n}$.

\emph{Remark 4.1:} Because the proximal operators $prox_{f^{1}_{i}}(\cdot)$ and $prox_{f^{2}_{i}}(\cdot)$, $i \in \lbrace 1, \cdots, n \rbrace$ are continuous and nonexpansive, the proposed algorithm is locally Lipschitz continuous even though $f^{1}_{i}(x_{i})$ and $f^{2}_{i}(x_{i})$ in problem ($\ref{Problem 2}$) are nonsmooth, which means the \textbf{smoothness} of the algorithm ($\ref{Algorithm 1}$) is guaranteed.

Different from the method in $\cite{CTOS3}$, algorithm ($\ref{Algorithm 1}$) is a fully distributed primal-dual algorithm to solve the saddle point dynamics of the Lagrangian function $L(x,v) = \sum_{j=0}^{2} \textbf{1}^{T}F^{j}(x) + v^{T}L_{nq}x + \frac{1}{2} x^{T} L_{nq} x$, where $F^{j}(x) = [f^{j}_{1}(x_{1}), f^{j}_{2}(x_{2}), \cdots, f^{j}_{n}(x_{n})]^{T}$, $j \in \lbrace 0, 1, 2 \rbrace$.

Algorithm ($\ref{Algorithm 1}$) can also be written in a compact form 
\begin{subequations}
\begin{align}
\dot{x}(t) = & PROX_{F^{1}}[x(t) - \nabla F^{0}(x(t)) \notag \\
& - \alpha L_{nq} v(t) - \alpha L_{nq} x(t) + \gamma z(t)] - x(t) \label{another proximal} \\
\dot{z}(t) = & PROX_{F^{2}}[x(t) - \gamma z(t)] - x(t) \label{estimator} \\
\dot{v}(t) = & \alpha L_{nq} x(t) 
\end{align}
\label{Algorithm 2}
\end{subequations}

\noindent where for any $\xi = [\xi_{1}^{T}, \xi_{2}^{T}, \cdots, \xi_{n}^{T}]^{T} \in \mathbb{R}^{nq}$, $\xi_{i} \in \mathbb{R}^{q}$, $i \in \lbrace 1, 2, \cdots, n \rbrace$, $j \in \lbrace 1, 2 \rbrace$, $PROX_{F^{j}}[\xi] = [(prox_{f^{j}_{1}}[\xi_{1}])^{T}, (prox_{f^{j}_{2}}[\xi_{2}])^{T}, \cdots, (prox_{f^{j}_{n}}[\xi_{n}])^{T}]^{T}$. $v = [v_{1}^{T}, v_{2}^{T}, \cdots, v_{n}^{T}]^{T}$ and $z = [z_{1}^{T}, z_{2}^{T}, \cdots, z_{n}^{T}]^{T}$.

\emph{Remark 4.2:} From (\ref{estimator}), it shows that $z(t)$ is a proximal-based estimator of $\partial F^{2}(x)$. With the help of estimator $z(t)$, another proximal operator (\ref{another proximal}), which employs the information of $z(t)$ instead of $\partial F^{2}(x)$, is presented to deal with the difficulty caused by the unproximable property of $F^{1}(x)+F^{2}(x)$. This estimation-based operator combined by (\ref{another proximal}) and (\ref{estimator})  is called the \textbf{double proximal operator}.

\emph{Lemma 4.1:}
Under the Assumptions 3.1-3.4, $(x^{*}, z^{*}, v^{*}) \in (\mathbb{R}^{nq}, \mathbb{R}^{nq}, \mathbb{R}^{nq})$ is an equilibrium of algorithm ($\ref{Algorithm 2}$) if and only if $x^{*}$ is a solution to problem ($\ref{Problem 2}$).

\begin{proof}
\textbf{(i)} Sufficiency: if $(x^{*}, z^{*}, v^{*})$ is an equilibrium of algorithm ($\ref{Algorithm 2}$), then according to the property ($\ref{Proximal Property}$) of proximal operator and algorithm ($\ref{Algorithm 2}$), one can have that
\begin{eqnarray}
- \nabla F^{0}(x^{*}) - \alpha L_{nq} v^{*} - \alpha L_{nq} x^{*} + \gamma z^{*} \!\!\!\!& \in \!\!\!\!& \partial F^{1}(x^{*}) \notag \\
- \gamma z^{*} \!\!\!\!& \in \!\!\!\!& \partial F^{2}(x^{*}) \notag \\
L_{nq} x^{*} \!\!\!\!& = \!\!\!\!& \textbf{0}_{nq}
\end{eqnarray}
which means $\textbf{0}_{nq} \in \nabla F^{0}(x^{*}) + \partial F^{1}(x^{*}) + \partial F^{2}(x^{*}) + L_{nq} v^{*}$.

According to Lemma 3.1, $x^{*}$ is a solution of problem ($\ref{Problem 2}$).

\textbf{(ii)} Necessity: if $x^{*}$ is a solution of the problem ($\ref{Problem 2}$), according to Lemma 3.1, there exist a point $z^{*} \in \mathbb{R}^{nq}$ and a constant $\gamma \in \mathbb{R}^{+}$ such that 
\begin{eqnarray}
\alpha L_{nq} x^{*} \!\!\!\!& = \!\!\!\!& \textbf{0}_{nq}, \notag\\
-\gamma z^{*} \!\!\!\!& \in \!\!\!\!& \partial F^{2}(x^{*}), \notag \\
- \nabla F^{0}(x^{*}) - L_{nq} v^{*} + \gamma z^{*} \!\!\!\!& \in \!\!\!\!& \partial F^{1}(x^{*}). \label{LN}
\end{eqnarray}

By combining the property ($\ref{Proximal Property}$) of the proximal operator with ($\ref{LN}$), we have
\begin{eqnarray}
PROX_{F_{2}}[x^{*} - \gamma z^{*}] - x^{*} \!\!\!\!\!\!\!\!& = \!\!\!\!\!\!& \textbf{0}_{nq}, \notag \\
PROX_{F^{1}}\![x^{*}\!\!\!\! -\!\! \nabla \! F^{0}\!(\!x^{*}\!)\!\! -\!\! \alpha L_{nq} v^{*}\!\!\!\! -\!\! \alpha L_{nq} x^{*}\!\! +\!\! \gamma z^{*}\!]\!\! -\!\! x^{*} \!\!\!\!\!\!\!\!& = \!\!\!\!\!\!& \textbf{0}_{nq}, 
\end{eqnarray}
which means that ($x^{*}, z^{*}, v^{*}$) is an equilibrium of algorithm ($\ref{Algorithm 2}$).
\end{proof}

\section{Main result}
In this section, we state the convergence result of the proposed distributed algorithm. Let $(x^{*}, z^{*}, v^{*})$ be an equilibrium of algorithm ($\ref{Algorithm 2}$). Define the Lyapunov candidate $V(x,z,v) = V_{1}(x,z,v) + V_{2}(x) + V_{3}(x,v)$, where

\begin{eqnarray}
V_{1}(x,z,v) \!\!\!\!& = \!\!\!\!& \frac{1}{2}\Vert x - x^{*} \Vert^{2} + \frac{1}{2} \gamma \Vert z - z^{*} \Vert^{2} \nonumber \\ 
\!\!\!\!& \!\!\!\!& - \gamma (x - x^{*})^{T}(z - z^{*}) + \frac{1}{2}\Vert v - v^{*}\Vert^{2} \notag \\
V_{2}(x) \!\!\!\!& = \!\!\!\!& \textbf{1}^{T}_{n}(F^{0}(x) - F^{0}(x^{*})) - (x - x^{*})^{T} \nabla F^{0}(x^{*}) \notag \\
\!\!\!\!& \!\!\!\!& + \frac{1}{2} \alpha x^{T} L_{nq} x \notag \\
V_{3}(x,v) \!\!\!\!& = \!\!\!\!& \alpha x^{T} L_{nq} (v - v^{*}) 
\end{eqnarray}

\emph{Theorem 5.1:} Consider algorithm ($\ref{Algorithm 1}$) or the equivalent form ($\ref{Algorithm 2}$). Suppose Assumptions 3.1-3.4 hold. 

\textbf{(i)} $V(x,z,v)$ is positive definite, radically unbounded, $V(x,z,v) = 0$ if and only if $(x,z,v) = (x^{*},z^{*},v^{*})$. 

\textbf{(ii)} Every equilibrium $(x^{*}, z^{*}, v^{*})$ of algorithm ($\ref{Algorithm 1}$) or ($\ref{Algorithm 2}$) is Lyapunov stable, and the trajectory of $(x(t),z(t),v(t))$ is bounded.

\textbf{(iii)} Moreover, the trajectory of $(x(t),z(t),v(t))$ converges and $\lim_{t \to \infty} x(t)$ is the solution to problem ($\ref{Problem 2}$).

\begin{proof}
\textbf{(i)} Since graph $\mathcal{G}$ is undirected, then $L^{T}_{nq} x = L_{nq} x$. And note that $L_{nq} x^{*} = 0_{nq}$ followed by ($\ref{Algorithm 2}$). It can be
easily verified that $V(x^{*}, z^{*},  v^{*}) = 0$.

Next, we show $V(x,z,v) > 0$ for all $(x,z,v) \neq (x^{*},z^{*},v^{*})$.

Since $f^{0}_{i}(x), i \in \lbrace 1,\cdots, n \rbrace$ are all convex, $\textbf{1}^{T}_{n}(F^{0}(x) - F^{0}(x^{*})) - (x - x^{*})^{T} \nabla F^{0}(x^{*}) \geq 0$. Note that $\lambda_{i}(L_{nq}) \geq 0, i \in \lbrace 1, \cdots, n \rbrace$. Hence 
\begin{eqnarray}
V_{2}(x) \!\!\!\!& = \!\!\!\!& \textbf{1}^{T}_{n}\!(F^{0}(x) \! - \! F^{0}(x^{*})) \!\! - \!\! (\!x \! - \! x^{*})^{T} \nabla F^{0}(\! x^{*}) + \frac{1}{2} \alpha x^{T} \!\! L_{nq} x \notag \\
\!\!\!\!& \geq \!\!\!\!& \textbf{1}^{T}_{nq}(F^{0}(x) - F^{0}(x^{*})) \!\! - \!\! (\!x \! - \! x^{*})^{T} \nabla F^{0}(\! x^{*}) \notag \\
\!\!\!\!& \geq \!\!\!\!& 0 \label{V2}
\end{eqnarray}

In addition, since $L_{nq}$ is positive semi-definite,
\begin{equation}
((x - x^{*}) + (v - v^{*}))^{T} L_{nq} ((x - x^{*}) + (v - v^{*})) \geq 0
\end{equation}

Therefore we have the conclusion that
\begin{eqnarray}
V_{3}(x,v) \!\!\!\!& = \!\!\!\!& \alpha x^{T} L_{nq} (v - v^{*}) \notag \\
\!\!\!\!& \geq \!\!\!\!& - \frac{\alpha}{2} (x  - x^{*})^{T} L_{nq} (x - x^{*}) \notag \\
\!\!\!\!& \!\!\!\!& - \frac{\alpha}{2} (v - v^{*})^{T} L_{nq} (v - v^{*}) \label{V31}
\end{eqnarray}

Let $\mu_{i}, i \in \lbrace 1, \cdots, n \rbrace$ be the eigenvalues of $L_{n} \in \mathbb{R}^{n\times n}$. Since the eigenvalues of $I_{q}$ are 1, it follows from the properties of Kronecker product that the eigenvalues of $L_{nq}$ are $\mu_{i} \times 1, i \in \lbrace 1, \cdots, n \rbrace$. Thus, $\lambda_{max}(L_{nq}) = \lambda_{max}(L_{n})$.

Note that $0 < \alpha < \frac{1}{\lambda_{max}(L_{n})}$ and $y^{T} L_{n} y \leq \lambda_{max}(L_{n}) \Vert y \Vert^{2}$. Hence 
\begin{eqnarray}
V_{3}(x,v) \!\!\!\!& \geq \!\!\!\!& -\frac{1}{2} \alpha \lambda_{max}(L_{n}) \Vert x - x^{*} \Vert^{2} \notag \\
\!\!\!\!& \!\!\!\!& - \frac{1}{2} \alpha \lambda_{max}(L_{n}) \Vert v - v^{*} \Vert^{2} \label{V32}
\end{eqnarray}

Since $0 < \gamma < 1 - \alpha \lambda_{max}(L_{n})$, therefore
\begin{eqnarray}
V_{1}(x,z,v) + V_{3}(x,v) \!\!\!\!& \geq \!\!\!\!& h_{1}\Vert x - x^{*} \Vert^{2} + \frac{1}{2} \gamma \Vert z - z^{*} \Vert^{2} \notag \\
\!\!\!\!& \!\!\!\!& - \gamma (x - x^{*})^{T}(z - z^{*}) \notag \\
\!\!\!\!& \!\!\!\!& + h_{1} \Vert v - v^{*} \Vert^{2} \notag \\
\!\!\!\!& = \!\!\!\!& h_{1} \left[ \Vert x - x^{*} \Vert^{2} + h_{2} \Vert z - z^{*} \Vert^{2} \right. \notag \\
\!\!\!\!& \!\!\!\!& \left. - 2h_{2} (x - x^{*})^{T}(z - z^{*}) \right] \notag \\
\!\!\!\!& \!\!\!\!& + h_{1} \Vert v - v^{*} \Vert^{2} \notag \\
\!\!\!\!& = \!\!\!\!& h_{1} \Vert (x - x^{*}) - h_{2}(z - z^{*})\Vert^{2} \notag \\
\!\!\!\!& \!\!\!\!& + h_{3} \Vert z - z^{*} \Vert^{2} + h_{1} \Vert v - v^{*} \Vert^{2} \notag \\
\!\!\!\!& \geq \!\!\!\!& 0 \label{V1+V3}
\end{eqnarray}
where $h_{1} = \frac{1}{2}(1 - \alpha \lambda_{max}(L_{n})) >0$, $h_{2} = \frac{\gamma}{1 - \alpha \lambda_{max}(L_{n})} >0$, and $h_{3} = \frac{1}{2} \gamma (1 - \frac{\gamma}{1 - \alpha \lambda_{max}(L_{n})}) >0$.

In the view of ($\ref{V2}$) and ($\ref{V1+V3}$), $V(x,z,v) \geq h_{1} \Vert (x - x^{*}) - h_{2}(z - z^{*})\Vert^{2} + h_{3} \Vert z - z^{*} \Vert^{2} + h_{1} \Vert v - v^{*} \Vert^{2} $. Clearly $V(x,z,v)$ is positive definite, radically unbounded, $V(x,z,v) \geq 0$ and is zero if and only if $(x,z,v) = (x^{*},z^{*},v^{*})$. 
\newline

\textbf{(ii)} It follows from algorithm ($\ref{Algorithm 2}$) that
\begin{eqnarray}
x + \dot{x} \!\!\!\!& = \!\!\!\!& PROX_{F^{1}}[x - \nabla F^{0}(x) - \alpha L_{nq} v - \alpha L_{nq} x + \gamma z] \notag \\
x^{*} \!\!\!\!& = \!\!\!\!& PROX_{F^{1}} [x^{*} - \nabla F^{0}(x^{*}) - \alpha L_{nq} v^{*} + \gamma z^{*}] \notag \\
x + \dot{z} \!\!\!\!& = \!\!\!\!& PROX_{F^{2}}[x - \gamma z] \notag \\
x^{*} \!\!\!\!& = \!\!\!\!& PROX_{F^{2}} [x^{*} - \gamma z^{*}] \label{x-x*}
\end{eqnarray}

Since $f^{1}_{i}(\cdot)$ and $f^{2}_{i}(\cdot)$ are both convex, $\partial f^{1}_{i}(\cdot)$ and $\partial f^{2}_{i}(\cdot)$ are both monotone for each agent $i$. According to the property ($\ref{Proximal Property}$) of the proximal operator, it follows from ($\ref{x-x*}$) that
\begin{subequations}
\begin{align}
& (\gamma (z - z^{*}) - (\nabla F^{0}(x) - \nabla F^{0}(x^{*})) - \alpha L_{nq} (v - v^{*}) \notag \\
& - \alpha L_{nq}(x - x^{*}) - \dot{x})^{T}((x - x^{*})) + \dot{x}) \geq 0\\
\notag \\
& (-\gamma (z - z^{*}) - \dot{z})^{T}((x - x^{*}) + \dot{z}) \geq 0
\end{align}
\label{Prox_pro}
\end{subequations}

From ($\ref{Prox_pro}$), we can deduce that
\begin{subequations}
\begin{align}
& \gamma (z - z^{*})^{T}(x - x^{*}) - (\nabla F^{0}(x) - \nabla F^{0}(x^{*}))^{T}(x - x^{*}) \notag \\
& - \alpha (v - v^{*})^{T} L_{nq} (x - x^{*})  - \alpha (x - x^{*})^{T} L_{nq} (x - x^{*}) \notag \\
& - \alpha (v - v^{*})^{T} L_{nq}  \dot{x} - \alpha (x - x^{*})^{T} L_{nq}  \dot{x} - (x \! - \! x^{*})^{T} \! \dot{x} \notag \\
& + \gamma (z \! - \! z^{*})^{T} \! \dot{x} \! - \! (\nabla F^{0}\! (x) \! - \! \nabla F^{0}\! (x^{*}\! ))^{T} \! \dot{x} \! - \! \Vert \dot{x} \Vert^{2} \geq 0\\
\notag \\
& -\gamma (z - z^{*})^{T} (x - x^{*}) - (x - x^{*})^{T} \dot{z} \notag \\
&- \gamma (z - z^{*})^{T} \dot{z} - \Vert \dot{z} \Vert^{2} \geq 0
\end{align}
\label{Proxmal relation}
\end{subequations}

The derivative of Lyapunov candidate $V(x,z,v)$ along the trajectory of algorithm ($\ref{Algorithm 1}$) or ($\ref{Algorithm 2}$) satisfies
\begin{eqnarray}
\dot{V}(x,z,v)\!\!\!\!&  = \!\!\!\!& (\! x \! - \! x^{*}\! )^{T} \! \dot{x} \! + \! \gamma (\! z \! - \! z^{*})^{T} \! \dot{z} \! - \! \gamma ((\! x \! - \! x^{*}\! )^{T} \! \dot{z} \! + \! (\! z \! - \! z^{*})^{T} \! \dot{x}) \notag \\
\!\!\!\!& \!\!\!\!& + \alpha (v \! - \! v^{*})^{T} L_{nq} (x \! - \! x^{*}) + (\nabla F^{0}(x) \! - \! \nabla F^{0}(x^{*}))^{T} \dot{x} \notag \\
\!\!\!\!& \!\!\!\!& + \alpha x^{T} \! L_{nq} \dot{x} \! + \! \alpha \dot{x}^{T} \! L_{nq} (v \! - \! v^{*}) \! + \! \alpha^{2} x^{T} \! L^2_{nq} x \label{dot_V}
\end{eqnarray}

According to (\ref{Proxmal relation}) and (\ref{dot_V}),
\begin{eqnarray}
\dot{V}(x,z,v) \!\!\!\!& \leq \!\!\!\!& - \Vert \dot{x} \Vert^{2} -\Vert \dot{z} \Vert^{2} - \alpha  x^{T} L_{M} x \notag \\
\!\!\!\!& \!\!\!\!& - (\nabla F^{0}(x) \! - \! \nabla F^{0}(x^{*}))^{T} \! (x\! -\! x^{*}) \notag \\
\!\!\!\!& \!\!\!\!& - (1 + \gamma)(x - x^{*})^{T} \dot{z}
\end{eqnarray}
where $L_{M} = L_{nq}^{T}(I_{nq} - \alpha L_{nq})$.

According to the Assumption 3.1, there exists a parameter $\beta > 0$ such that
\begin{eqnarray}
-(1 + \gamma)(x - x^{*})\dot{z} \!\!\!\!& \leq \!\!\!\!& \frac{1}{2}(1\! +\! \gamma) \frac{1}{\beta}(x - x^{*})^{2} + \frac{1}{2}(1 + \gamma) \beta \dot{z}^{2} \notag \\
\!\!\!\!& \leq \!\!\!\!& \frac{1}{2}\!(1\!\! +\!\! \gamma) \frac{1}{c\beta}(\nabla \! F^{0}(x)\!\! -\!\! \nabla \! F^{0}(x^{*}))\!^{T}\!(\! x\!\! -\!\! x^{*}\!) \notag \\
\!\!\!\!& \!\!\!\!& + \frac{1}{2}(1 + \gamma) \beta \dot{z}^{2} 
\end{eqnarray}

Then we have the conclusion that
\begin{eqnarray}
\dot{V}(x,z,v) \!\!\!\!& \leq \!\!\!\!& -\Vert \dot{x} \Vert^{2} - b_{1} \Vert \dot{z} \Vert^{2} - \alpha x^{T} L_{M} x \notag \\
\!\!\!\!& \!\!\!\!& - b_{2} (\nabla F^{0}(x) - \nabla F^{0}(x^{*}))^{T}(x - x^{*})
\end{eqnarray}
where $b_{1} = 1 - \frac{1}{2}(1 + \gamma) \beta$, $b_{2} = 1 - \frac{1}{2}(1 + \gamma) \frac{1}{c\beta}$.

In order to illustrate that there exists $\beta > 0$ to make $b_{1} > 0$ and $b_{2} > 0$, here we define a function $G(\gamma)$ of $\gamma$ and its derivative $\frac{d G}{d \gamma}$ as follows
\begin{eqnarray}
G(\gamma) \!\!\!\!& = \!\!\!\!& \frac{2}{\gamma + 1} - \frac{\gamma + 1}{2c} \notag \\
\frac{d G(\gamma)}{d \gamma} \!\!\!\!& = \!\!\!\!& -\frac{2}{(\gamma + 1)^{2}} - \frac{1}{2c} < 0 \label{Beta}
\end{eqnarray}

Note that $0 < \gamma < 1 - \alpha \lambda_{max}(L_{n}) < 1$ and $c >1$. According to  ($\ref{Beta}$), we have
\begin{equation} 
G_{min}(\gamma) > G(1) = 1 - \frac{1}{c} > 0 
\end{equation}

As the result, there exists $\beta$ such that $\frac{1 + \gamma}{2c} < \beta < \frac{2}{1 + \gamma}$, which means that
\begin{eqnarray}
b_{1} \!\!\!\!& = \!\!\!\!& 1 - \frac{1}{2}(1 + \gamma) \beta > 0 \notag \\
b_{2} \!\!\!\!& = \!\!\!\!& 1 - \frac{1}{2}(1 + \gamma) \frac{1}{c\beta} > 0
\end{eqnarray}

Since $f_{i}^{0}, i \in \lbrace 1, \cdots, n \rbrace$ are all convex functions, which means $(\nabla F^{0}(x) - \nabla F^{0}(x^{*}))^{T}(x - x^{*}) \geq 0$, as the result,
\begin{equation}
\begin{array}{l}
\dot{V}(x,z,v) \leq -\Vert \dot{x} \Vert^{2} - b_{1} \Vert \dot{z} \Vert^{2} - \alpha x^{T} L_{M} x\\
\end{array}
\end{equation}

Note that $\lambda_{i}(L_{M})= \lambda_{i}(L_{nq} (I_{nq} - \alpha L_{nq})) \geq 0$, $i \in \lbrace 1, \cdots, n \rbrace$. Hence
\begin{equation}
\dot{V}(x,z,v) \leq -\Vert \dot{x} \Vert^{2} - b_{1} \Vert \dot{z} \Vert^{2} \leq 0
\end{equation}

Additionally, since $V(x,z,v)$ is positive-definite, radically unbounded, lower bounded, $(x^{*},z^{*}, v^{*})$ is Lyapunov stable and the trajectory $(x(t),z(t),v(t))$ is bounded.

\textbf{(iii)} Define
\begin{equation}
\begin{array}{ll}
R &\ \!\!\!\!\!\! = \lbrace (x,z,v) : \dot{V}(x,z,v) = 0 \rbrace\\
&\ \!\!\!\!\!\! \subseteq \lbrace (x,z,v) : \dot{x} = \textbf{0}_{nq}, \dot{z} = \textbf{0}_{nq} \rbrace\\
\end{array}
\end{equation}

Let $M$ be the largest invariant set of $R$. It follows from the invariance principle that $(x(t),z(t),v(t)) \to M$ as $t \to \infty$. Assume $(\bar{x}(t), \bar{z}(t), \bar{v}(t))$ is a trajectory of algorithm ($\ref{Algorithm 1}$) or ($\ref{Algorithm 2}$) such that $(\bar{x}(0), \bar{z}(0), \bar{v}(0)) \in M$. Then $(\bar{x}(t), \bar{z}(t), \bar{v}(t)) \subset M$ for all $t \geq 0$. Therefore, $\dot{\bar{x}}(t) \equiv \textbf{0}_{nq}$, $\dot{\bar{z}}(t) \equiv \textbf{0}_{nq}$ and
\begin{equation}
\dot{\bar{v}}(t) \equiv L\bar{x}(0) \equiv C
\end{equation}

If $C \neq 0$, $\bar{v}(t)$ will be unbounded, which is a contradiction. Hence $M \subseteq \lbrace (x,z,v) : \dot{x} = \textbf{0}_{nq}, \dot{z} = \textbf{0}_{nq}, \dot{v} = \textbf{0}_{nq} \rbrace$. By Part (ii), every point in $M$ is Lyapunov stable. It follows from Lemma 3.1 that $(x(t), z(t), v(t))$ converges to an equilibrium point. Due to Lemma 4.1, $\lim_{t \to \infty} x(t)$ is a solution to problem (\ref{Problem 2}).
\end{proof}

Theorem 5.1 shows that every equilibrium of algorithm ($\ref{Algorithm 1}$) or ($\ref{Algorithm 2}$) is Lyapunov stable and every state trajectory converges to one equilibrium point, which is the optimal point of problem ($\ref{Problem 2}$). Moreover, by using the double proximal operator, algorithm ($\ref{Algorithm 1}$) or ($\ref{Algorithm 2}$) is a smooth algorithm and its analysis is easier compared with the nonsmooth algorithm in $\cite{DCO4}$. This method gives new ideas for designing distributed algorithms for nonsmooth optimizations without using nonsmooth analysis.

\section{Simulations}
In this section, simulations are performed to validate our proposed distributed optimization algorithm.

Consider the problem ($\ref{Problem 2}$) with four agents moving in a 2-D space with first-order dynamics,
where $x_{i} = [ x^{1}_{i}, x^{2}_{i} ]^{T} \in \mathbb{R}^{2}, i \in \lbrace 1, 2, 3, 4 \rbrace$. The local cost function $f_{i}(x_{i})$ for agent $i$ consists of the follows
\begin{equation}
\begin{array}{l}
f^{0}_{i}(x_{i}) = \Vert x_{i} - m_{i}\Vert^{2}\\
f^{1}_{i}(x_{i}) =
\begin{cases} 
0,  & \mbox{if } x_{i} \in \Omega_{i}\\
\infty, & \mbox{if } x_{i} \notin \Omega_{i}
\end{cases}\\
f^{2}_{i}(x_{i}) = \Vert x_{i} - p_{i} \Vert_{1}\\
\end{array}
\label{Simulation_f}
\end{equation}
where $m_{i} = [m^{1}_{i}, m^{2}_{i}]^{T} = [i - 2.5, 0]^{T}$, $\Omega_{i} = \lbrace \delta \in \mathbb{R}^{2} \vert \Vert \delta - x_{i}(0) \Vert^{2} \leq 64\rbrace$, $p_{i} = [p^{1}_{i}, p^{2}_{i}]^{T} = [0, i - 2.5]^{T}$. The $f^{0}_{i}(x_{i})$, $f^{1}_{i}(x_{i})$ and $f^{2}_{i}(x_{i})$ represent respectively the quadratic objective, the indicator function of the constraint set $x_{i} \in \Omega_{i}$, and the $l_{1}$ penalty with an anchor $p_{i}$ for each agent $i$.

Based on ($\ref{Simulation_f}$), the gradient of $f^{0}_{i}$ and the proximal operators of $f^{1}_{i}$ and $f^{2}_{i}$ for agent $i$ are
\begin{eqnarray}
\nabla f^{0}_{i}(x_{i}) \!\!\!\!& = \!\!\!\!& [ 2(x^{1}_{i} - m^{1}_{i}), 2(x^{2}_{i} - m^{2}_{i}) ]^{T} \notag \\
prox_{f^{1}_{i}}[\eta_{1}] \!\!\!\!& = \!\!\!\!& \arg \min_{\delta \in \Omega_{i} } \Vert \delta - \eta_{1} \Vert^{2} \notag \\
prox_{f^{2}_{i}}[\eta_{2}] \!\!\!\!& = \!\!\!\!& [\phi^{1}_{i}, \phi^{2}_{i}]^{T} \notag \\
\end{eqnarray}
where $\eta_{1} \in \mathbb{R}^{2}$, $\eta_{2} \in \mathbb{R}^{2}$ and for $j \in \lbrace 1, 2 \rbrace$,
\begin{eqnarray}
\phi^{j}_{i} \!\!\!\!& = \!\!\!\!&
\begin{cases} 
\eta_{2}^{j}-1,  & \mbox{if } \eta_{2}^{j} > p^{j}_{i}+1 \\
\eta_{2}^{j},  & \mbox{if } p^{j}_{i}-1 \leq \eta_{2}^{j} \leq  p^{j}_{i}+1 \\
\eta_{2}^{j}+1,  & \mbox{if } \eta^{j}_{2} <  p^{j}_{i}-1
\end{cases} 
\end{eqnarray}

Note that the proximal operator of $f^{1}_{i}(x_{i})+f^{2}_{i}(x_{i}), i \in \lbrace 1, 2, 3, 4 \rbrace$ 
\begin{equation}
prox_{(f^{1}_{i}+f^{2}_{i})}[\eta_{3}]= \arg \min_{\delta \in \Omega_{i}} \lbrace \Vert \delta \! - \! p_{i} \Vert_{1} + \frac{1}{2} \Vert \delta - \eta_{3} \Vert^{2} \rbrace
\end{equation}
is not proximable, where $\eta_{3} \in \mathbb{R}^{2}$, hence the proximal algorithms $\cite{DCFO1}$-$\cite{DCFO3}$ may not fit for this problem.

The Laplacian of the undirected graph $\mathcal{G}$ is given by
\begin{equation}
L_{4} = 
\begin{bmatrix}
1      & -1 & 0   & 0 \\
-1      & 2 & -1   & 0 \\
0      & -1 & 2   & -1 \\
0      & 0 & -1   & 1 \\
\end{bmatrix}
\end{equation}
and the Laplacian spectrum of the topology is $\lbrace 0, 0.586, 2, 3.414 \rbrace$. According to the knowledge of $\lambda_{max}(L_{4})$, we set $\alpha = 0.2$ and $\gamma = 0.3$ as the coefficients in algorithm ($\ref{Algorithm 2}$). The initial positions of the agents 1, 2, 3, and 4 are set as $x_{1}(0) = [-4, 5.5]^{T}$, $x_{2}(0) = [6, 5]^{T}$, $x_{3}(0) = [5, -3.5]^{T}$, and $x_{4}(0) = [-5, -5]^{T}$. We set the initial values for the Lagrangian multipliers $v_{i}, i \in \lbrace 1, 2, 3, 4 \rbrace$ and auxiliary variables $z_{i}, i \in \lbrace 1, 2, 3, 4 \rbrace$ as zeros. The optimal solution to this problem is $x^{*} = [0,0]^{T}$.

\begin{figure}
\centering
\includegraphics[width=0.35\textwidth]{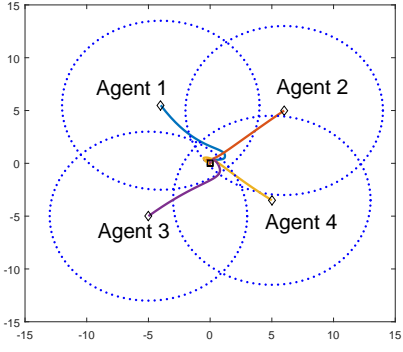}
\caption{The motions of the entire system in a 2-D space with algorithm ($\ref{Algorithm 2}$)}
\label{Fig.1}
\end{figure}

\begin{figure}
\centering
\subfigure{
\includegraphics[width=0.35\textwidth]{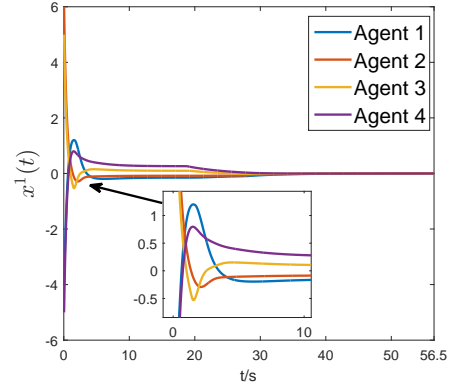}}
\subfigure{
\includegraphics[width=0.35\textwidth]{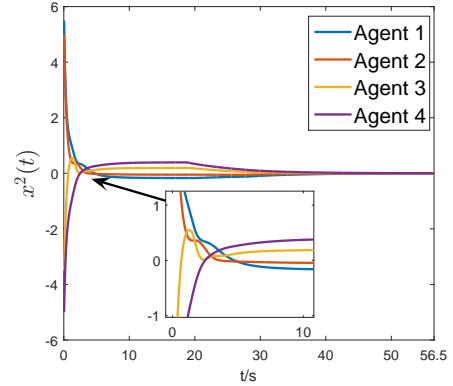}}
\caption{The trajectories of $x_{i}(t), i \in \lbrace 1, 2, 3, 4 \rbrace$ with algorithm ($\ref{Algorithm 2}$)}
\label{Fig.2}
\end{figure}

\begin{figure}[h]
\centering
\includegraphics[width=0.35\textwidth]{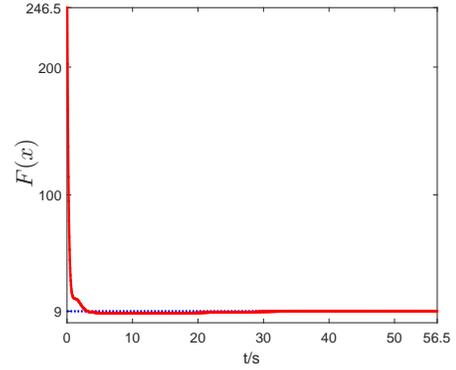}
\caption{The trajectory of $F(x)$ with algorithm ($\ref{Algorithm 2}$)}
\label{Fig.3}
\end{figure}

\begin{figure}
\centering
\subfigure{
\includegraphics[width=0.35\textwidth]{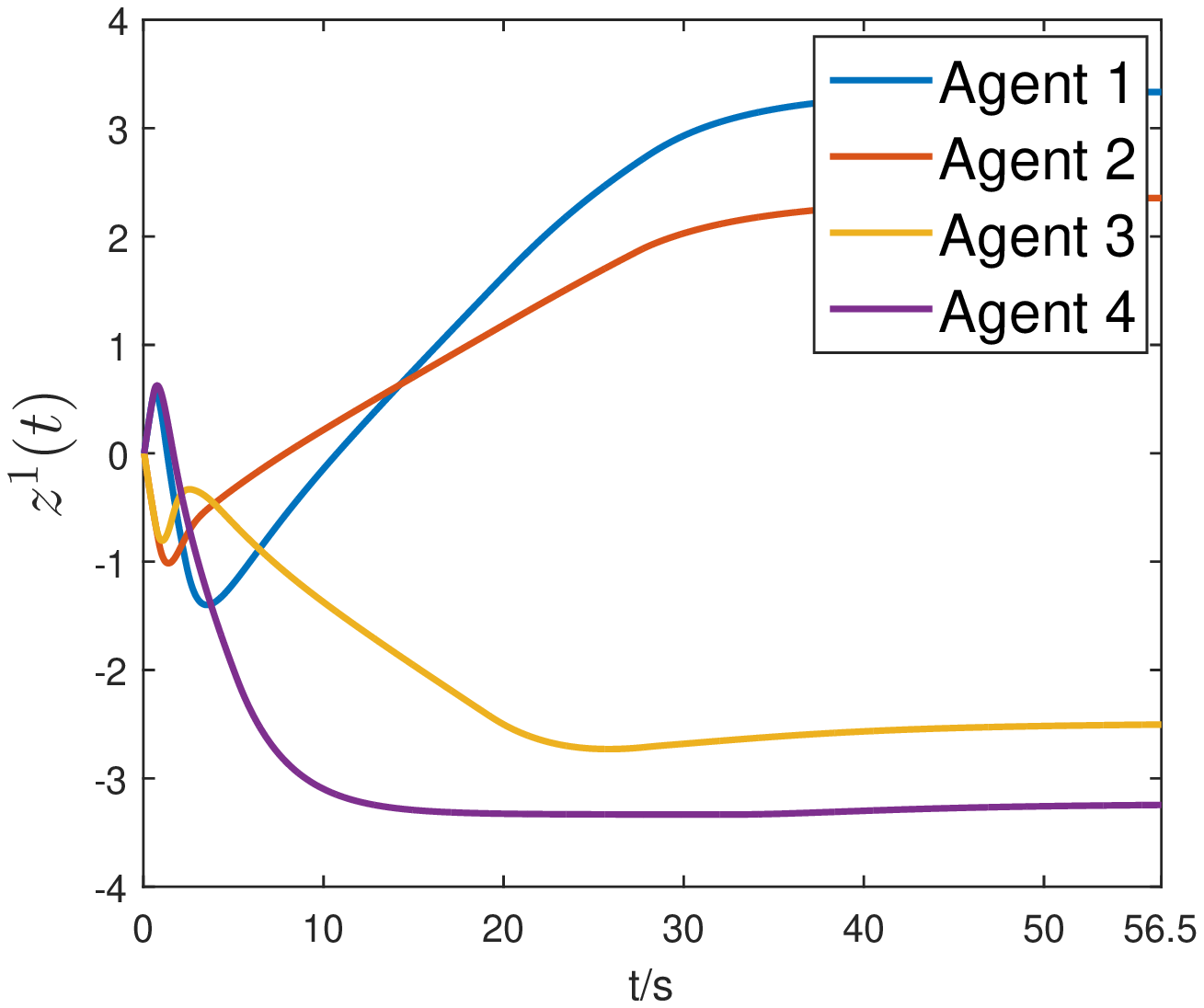}}
\subfigure{
\includegraphics[width=0.35\textwidth]{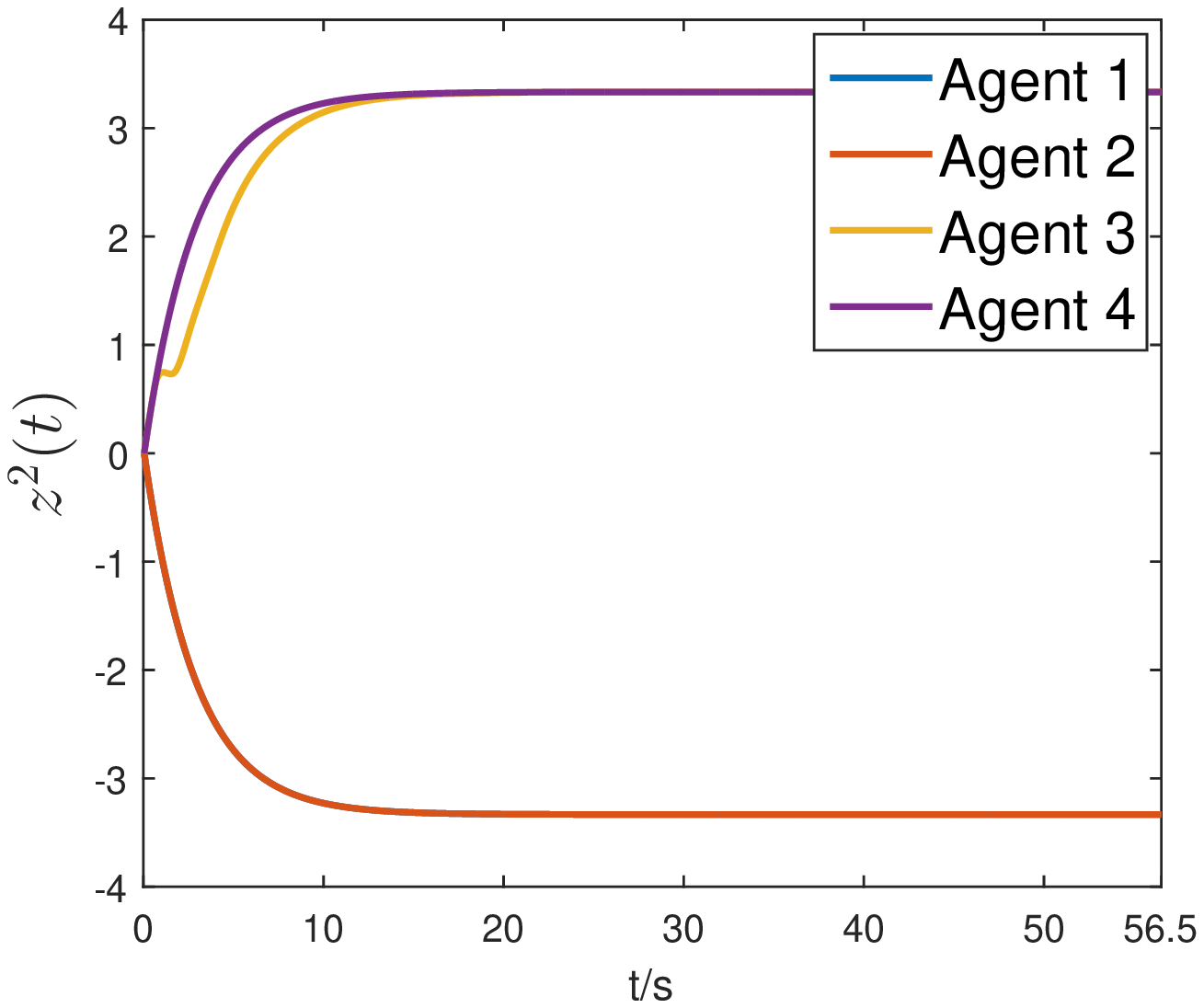}}
\caption{The trajectories of $z_{i}(t), i \in \lbrace 1, 2, 3, 4 \rbrace$ with algorithm ($\ref{Algorithm 2}$)}
\label{Fig.4}
\end{figure}

\begin{figure}
\centering
\subfigure{
\includegraphics[width=0.35\textwidth]{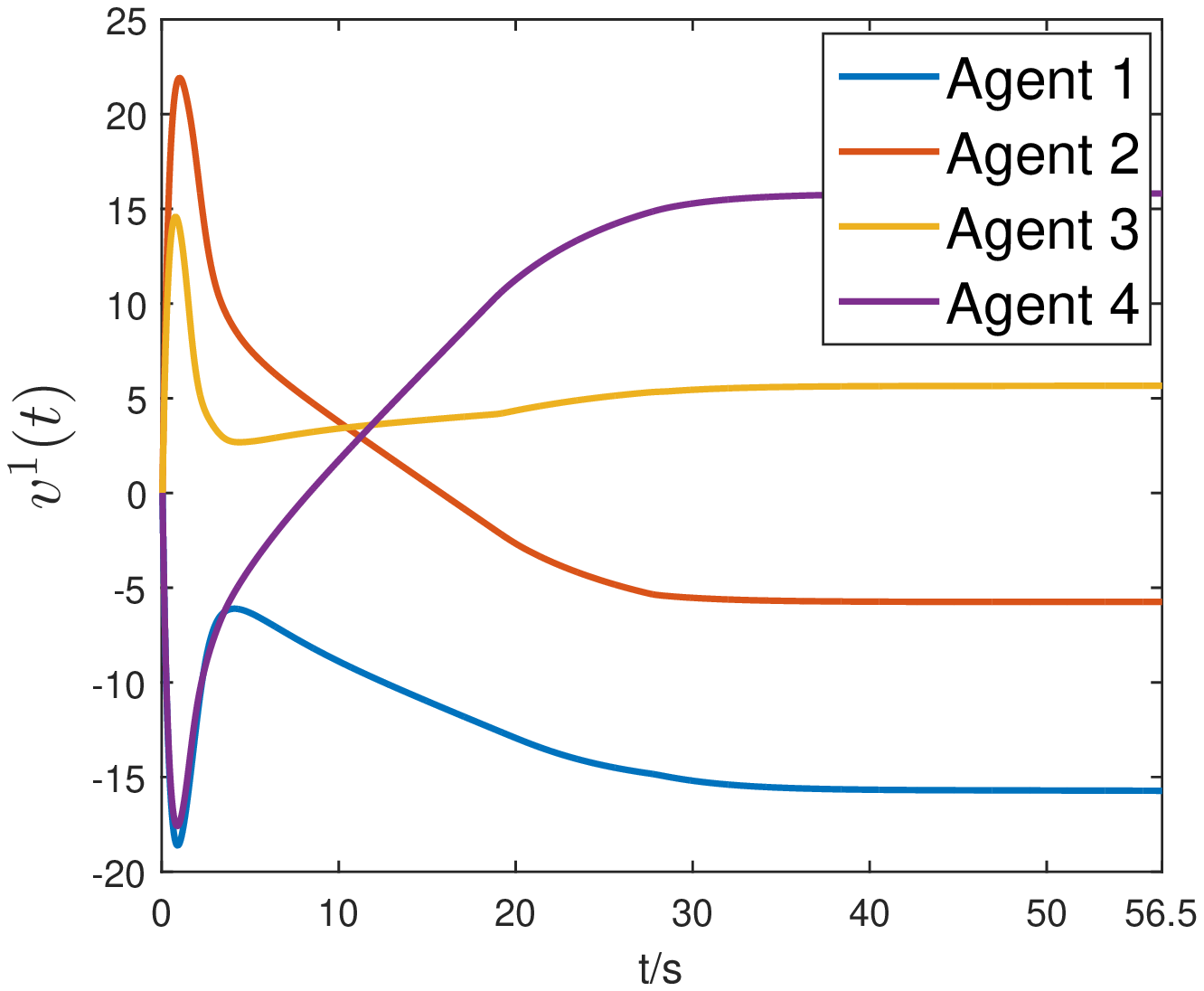}}
\subfigure{
\includegraphics[width=0.35\textwidth]{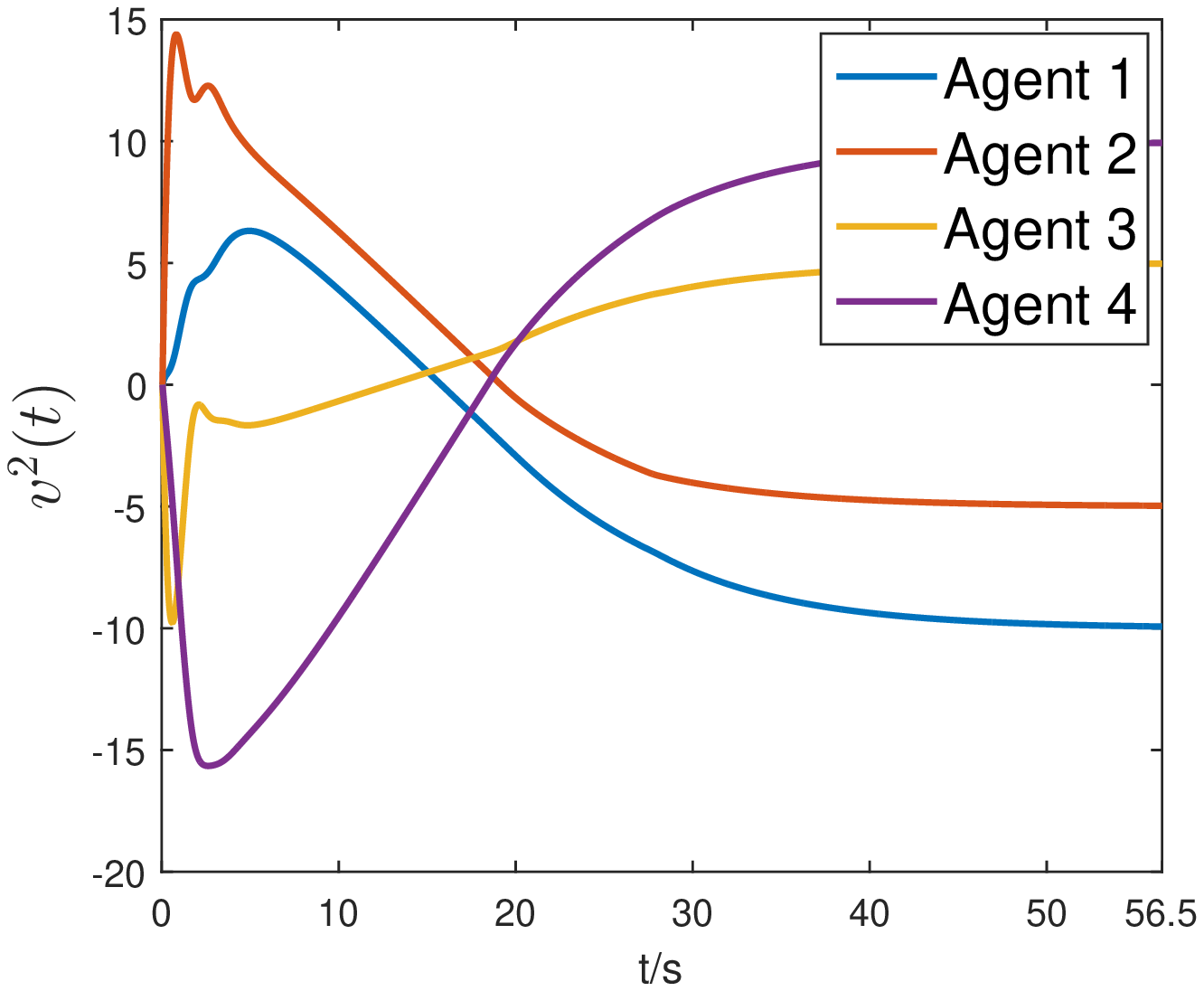}}
\caption{The trajectories of $v_{i}(t), i \in \lbrace 1, 2, 3, 4 \rbrace$ with algorithm ($\ref{Algorithm 2}$)}
\label{Fig.5}
\end{figure}

The motions of the entire system versus time are shown in Fig.$\ref{Fig.1}$. Fig.$\ref{Fig.2}$ gives the trajectories of $x_{i}(t), i \in \lbrace 1, 2, 3, 4 \rbrace$. Fig.$\ref{Fig.3}$ shows the trajectory of $F(x)$, which proves that the summation of the local cost functions is minimized. It can be seen from Fig.$\ref{Fig.1}$-Fig.$\ref{Fig.3}$ that all the agents converge to the same optimal solution which minimizes the summation of local cost functions. Fig.$\ref{Fig.4}$ and Fig.$\ref{Fig.5}$ show the trajectories of the auxiliary variables $z_{i}(t), i \in \lbrace 1, 2, 3, 4 \rbrace$ and the Lagrange multipliers $v_{i}(t), i \in \lbrace 1, 2, 3, 4 \rbrace$ respectively, which also verify the boundedness of algorithm ($\ref{Algorithm 2}$)'s trajectories. 

\begin{figure}
\centering
\subfigure{
\includegraphics[width=0.35\textwidth]{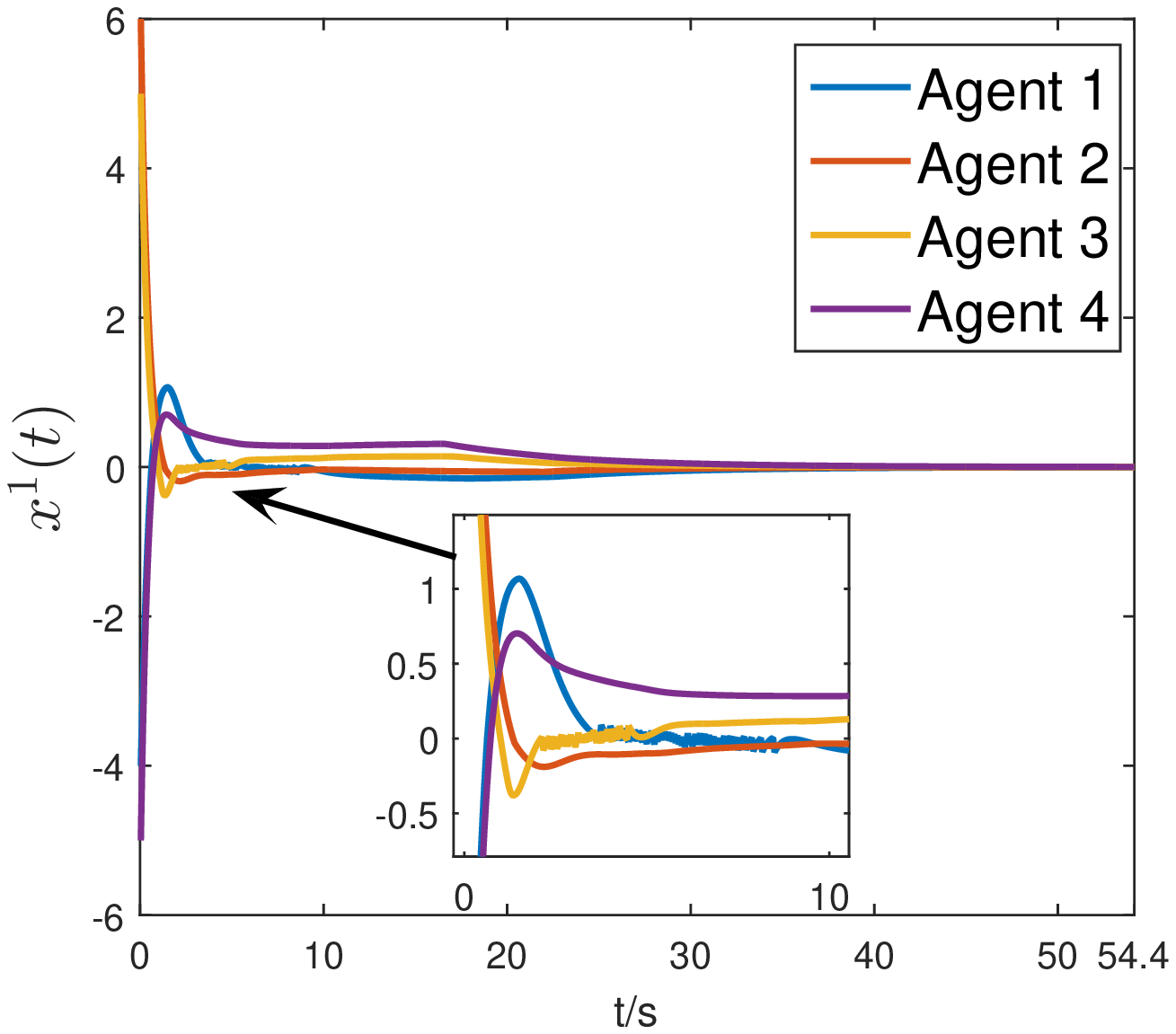}}
\subfigure{
\includegraphics[width=0.35\textwidth]{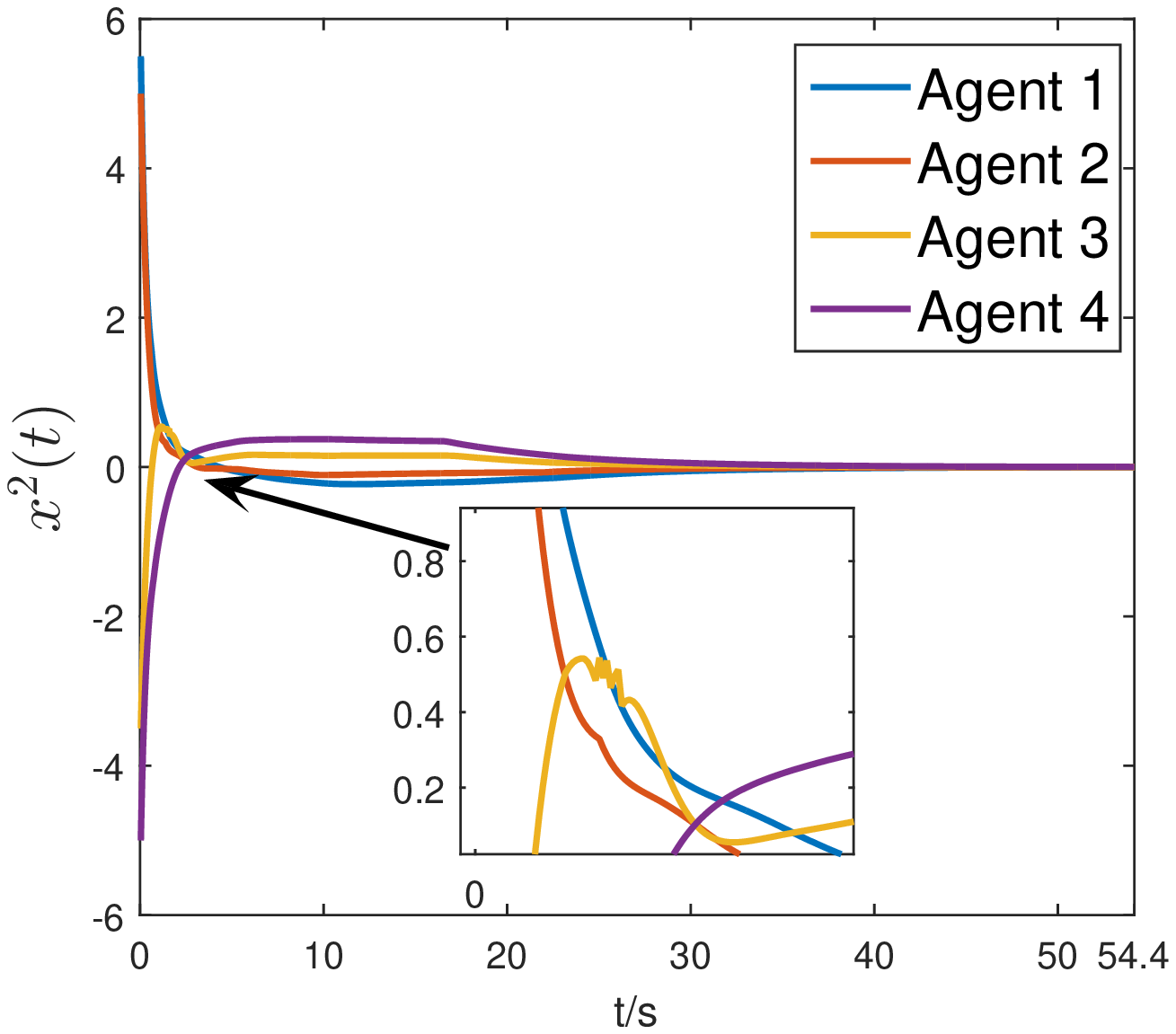}}
\caption{The trajectories of $x_{i}(t), i \in \lbrace 1, 2, 3, 4 \rbrace$ with the nonsmooth algorithm}
\label{Fig.6}
\end{figure}

Comparatively, Fig.$\ref{Fig.6}$ shows the trajectories of $x_{i}(t), i \in \lbrace 1, 2, 3, 4 \rbrace$ if the subgradient of $f^{2}_{i}(x_{i})$ is directly employed like the nonsmooth algorithm $\cite{DCO1}$ does. From Fig.$\ref{Fig.6}$, we can find that the trajectories of $x_{i}(t), i \in \lbrace 1, 2, 3, 4 \rbrace$ go zigzag from 0s to 10s, 
which is unacceptable in the distributed optimization problem implemented by physical systems.

\section{Conclusions}
In this note, a distributed smooth double proximal primal-dual algorithm is proposed to solve a class of distributed nonsmooth consensus optimization problem, which is called as \emph{single-smooth plus double-nonsmooth} (SSDN) problem. The double proximal operator which contains a proximal-based estimator is employed to tackle the difficulty caused by the unproximable property of the local cost functions and ensure smoothness of the proposed algorithm. Moreover, the proposed algorithm can also give a new viewpoint to many widely studied distributed constrained optimization problems. Future extensions will involve considering the SSDN problem with the directed and switching topologies.


\begin{thebibliography}{52}
\addtolength{\itemsep}{-1.5 em} 
\setlength{\itemsep}{0pt}
\footnotesize








\bibitem{CTOS3} E. K. Ryu and W. Yin, ``Proximal-proximal-gradient method," \emph{arXiv preprint arXiv:1708.06908}, 2017.

\bibitem{Lasso} R. Tibshirani, M. Saunders, S. Rosset, J. Zhu and K. Knight, ``Sparsity and smoothness via the fused LASSO", \emph{J. R. Stat. Soc. Ser. B-Stat. Methodol.}, vol. 67, no.1, pp. 91-108 , 2005.





\bibitem{DCOR1} P. Lin, W. Ren, and J. A. Farrell, ``Distributed continuous-time optimization: nonuniform gradient gains, finite-time convergence, and convex constraint set," \emph{IEEE Trans. Autom. Control}, vol. 62, no. 5, pp. 2239-2253, 2017.

\bibitem{Yipeng1} P. Yi, Y. Hong, and F. Liu, ``Distributed gradient algorithm for constrained optimization with application to load sharing in power systems," \emph{Syst. Control Lett.}, vol. 83, pp. 45-52, 2015.

\bibitem{DCON2} T. H. Chang, A. Nedi$\acute{c}$, and A. Scaglione, ``Distributed constrained optimization by consensus-based primal-dual perturbation method," \emph{IEEE Trans. Autom. Control}, vol. 59, no. 6, pp. 1524-1538, 2014.

\bibitem{Yipeng2} P. Yi, Y. Hong, and F. Liu, ``Initialization-free distributed algorithms for optimal resource allocation with feasibility constraints and application to economic dispatch of power systems," \emph{Automatica}, vol. 74, pp. 259-269, 2016.

\bibitem{DCO1} A. Nedi$\acute{c}$ and A. Ozdaglar, ``Distributed subgradient methods for multi-agent optimization", \emph{IEEE Trans. Autom. Control}, vol. 54, no. 1, pp. 48-61, 2009.


\bibitem{DCO111} D. Yuan, S. Xu, and H. Zhao, ``Distributed primal-dual subgradient method for multi-agent optimization via consensus algorithms," \emph{IEEE Trans. Syst., Man, Cyber., Cyber.}, vol. 41, no. 6, pp. 1715-1724, 2011.

\bibitem{DCO3} S. Lee and M. M. Zavlanos, ``Approximate projection methods for decentralized optimization with functional constraints," \emph{IEEE Trans. Autom. Control, doi:10.1109/TAC.2017.2778696}, 2017.

\bibitem{Xie1} Z. Qiu, S. Liu, and L. Xie, ``Distributed constrained optimal consensus of multi-agent systems," \emph{Automatica}, vol. 68, pp. 209-215, 2016.

\bibitem{Yang1} S. Yang, Q. Liu, and J. Wang, ``A multi-agent system with a proportional-integral protocol for distributed constrained optimization," \emph{IEEE Trans. Autom. Control}, vol. 62, no. 7, pp. 3461-3467, 2017.

\bibitem{DOC1} D. Mateos-N$\acute{u} \tilde{n}$ez, and J. Cort$\acute{e}$s, ``Distributed saddle-point subgradient algorithms with Laplacian averaging," \emph{IEEE Trans. Autom. Control}, vol. 62, no. 6, pp. 2720-2735, 2017.

\bibitem{DCO2} Q. Liu and J. Wang, ``A second-order multi-agent network for bound constrained distributed optimization," \emph{IEEE Trans. Autom. Control}, vol. 60, no. 12, pp. 3310-3315, 2015.

\bibitem{DCO4} X. Zeng, P. Yi, and Y. Hong, ``Distributed continuous-time algorithm for constrained convex optimizations via nonsmooth analysis approach," \emph{IEEE Trans. Autom. Control}, vol. 62, no. 10, pp. 5227-5233, 2017.

\bibitem{DCO5} S. Liang, X. Zeng, and Y. Hong, ``Distributed nonsmooth optimization with coupled inequality constraints via modified Lagrangian function," \emph{IEEE Trans. Autom. Control}, vol. 63, no.6, pp. 1753-1759, 2018.

\bibitem{NCO} Z. Denkowski, S. Mig$\acute{o}$rski, and N. S. Papageorgiou, \emph{An Introduction to Nonlinear Analysis: Theory}, New York, NY: Springer-Verlag, 2003.

\bibitem{DCFO1} W. Shi, Q. Ling, G. Wu, and W. Yin, ``A proximal gradient algorithm for decentralized composite optimization," \emph{IEEE Trans. Signal Process.}, vol. 63, no. 22, pp. 6013-6023, 2015.

\bibitem{DCFO4} Z. Li, W. Shi, and M. Yan, ``A decentralized proximal-gradient method with network independent step-sizes and separated convergence rates," \emph{arXiv preprint arXiv:1704.07807}, 2017.




\bibitem{DCFO2} M. Hong and T. H. Chang, ``Stochastic proximal gradient consensus over random network," \emph{IEEE Trans. Signal Process.}, vol. 65, no. 11, pp. 2933-2948, 2017.

\bibitem{DCFO3} N. S. Aybat, Z. Wang, T. Lin, and S. Ma, ``Distributed linearized alternating direction method of multipliers for composite convex consensus optimization," \emph{IEEE Trans. Autom. Control}, vol. 63, no. 1, pp. 5-20, 2018.


\bibitem{PM} N. Parikh and S. Boyd, ``Proximal algorithms," \emph{Foundations and Trends in Optimization}, vol. 1, no. 3, pp. 123-231, 2014.

\bibitem{TOS1} D. Davis and W. Yin, ``A three-operator splitting scheme and its optimization applications," \emph{Set-Valued Var. Anal.}, vol. 25, no. 4, pp. 829-858 , 2017.

\bibitem{TOS2} M. Yan, ``A new primal-dual algorithm for minimizing the sum of three functions with a linear operator," \emph{J. Sci. Comput.}, vol. 76, no. 3, pp. 1698-1717, 2018.

\bibitem{TOS3} F. Pedregosa and G. Gidel, ``Adaptive three operator splitting," \emph{arXiv preprint arXiv:1804.02339}, 2018.

\bibitem{NO} A. Ruszczynski, \emph{Nonlinear Optimization}, Princeton, NJ: Princeton University Press, 2006.








\end{thebibliography}
\end{document}